\documentclass[12pt]{amsart}

\usepackage{amsmath,amssymb,amscd,amsfonts,verbatim}
\newtheorem{thm}{Theorem}[section]

\newtheorem{prop}[thm]{Proposition}

\newtheorem{assu-nota}[thm]{Assumption--Notation}
\theoremstyle{remark}

\newcommand{\N}{\mathbb N}
\newcommand{\pp}{\mathbb P}

\newcommand{\epsi}{\varepsilon}
\newcommand{\al}{\alpha}
\newcommand{\be}{\beta}

\newcommand{\la}{\lambda}

\numberwithin{equation}{section}

\begin{document}
\title[The Severi inequality
$K^2\ge 4\chi$]{The Severi inequality\\
$K^2\ge 4\chi$\\for surfaces of maximal Albanese dimension}

\author{Rita Pardini}
\date{}
\begin{abstract}
We prove the so-called Severi inequality, stating  that the invariants of  a minimal smooth
complex  projective surface of maximal Albanese dimension satisfy: $$K^2_S\ge
4\chi(S).$$
\end{abstract}

\maketitle
\section{introduction}

In this paper we prove the so-called Severi inequality:
\smallskip  

\noindent{\em If $S$ is a minimal
smooth complex projective surface of maximal Albanese dimension, then $K^2_S\ge 4\chi(S)$}.
\medskip

This inequality has a long history. Severi (\cite{sev}) stated it  as a theorem  
in 1932, but his proof was not correct (cf. \cite{cat}). Around  the end of the   1970's the
inequality was posed  as a conjecture by  Reid (\cite{miles}) and by Catanese (\cite{cat}). Motivated
by this conjecture, Xiao wrote the foundational paper \cite{xiao} on  surfaces fibred over a curve,
in which he proves the conjecture in the special case of a surface admitting a fibration over
a curve of positive genus.  Finally,  at the end of the 1990's, the conjecture was almost solved by
Manetti (\cite{manetti}),  who proved  it  under  the additional assumption that
the surface have ample canonical bundle. His proof is based on a
very fine  study of the positivity properties of the tautological line bundle of
$\pp(\Omega^1_S\otimes\omega_S)$ and it allows him also to classify the surfaces for which equality
holds, still under the assumption that the canonical bundle  be ample.
\medskip

Our proof of the conjecture in the general case  is completely different from
Manetti's. We derive the Severi inequality from a well known inequality for the slope of a fibered
surface proven by Xiao in the same paper
\cite{xiao} and, independently,  by Cornalba--Harris (\cite{ch}) in the semistable case. In order to do this,   we  construct first an infinite sequence of fibred 
surfaces
$f_d\colon Y_d\to \pp^1$  such  that the slope of $f_d$  
converges to the ratio
$K^2_S/\chi(S)$ as $d$ goes to infinity (cf. Prop. \ref{trick}). Then we obtain the Severi 
inequality by applying the   slope inequality to the fibrations 
$f_d$ and taking the limit for $d\to \infty$.
\medskip 

\noindent{\bf Acknowledgments.} This work was partially supported by  P.R.I.N. 2002
``Geometria sulle variet\`a algebriche'' of Italian M.I.U.R.
I am indebted to Xiao Gang and J\'anos Koll\'ar for suggesting  simplifications of the original proof.
\medskip

\noindent{\bf Notation and conventions.}
We work over the field of complex numbers; all varieties are projective.
 Linear equivalence is denoted by $\equiv$ and numerical equivalence by $\sim_{num}$. We say that a
surface $S$ has maximal Albanese dimension if the image of its Albanese map is a surface.

The remaining 
notation is standard in algebraic geometry. We just recall the numerical invariants associated to a
surface
$S$: the self intersection $K^2_S$ of the canonical divisor $K_S$,  the {\em geometric genus}
$p_g(S):=h^0(S,K_S)$, the {\em irregularity}
$q(S):=h^0(S,\Omega^1_S)$ and the {\em Euler--Poincar\'e characteristic} $\chi(S):=p_g(S)-q(S)+1$.
\section{The Severi inequality}

In this section we  prove the following:
\begin{thm}[Severi inequality]\label{severi}
Let $S$ be a smooth minimal complex surface of maximal Albanese dimension. Then:
$$K^2_S\ge 4\chi(S).$$
\end{thm}
\bigskip

We begin by recalling some facts about surfaces fibred over a curve.

Let $Y$ be a smooth  complex surface and let $f\colon Y\to B$ be a fibration onto a smooth
curve of genus $b$. Assume that $f$ is relatively minimal, i.e., that no fibre of $f$ contains a
$-1-$curve, and that  a general fibre of $f$ has genus $g\ge 2$. It is well known (cf.
\cite{appendix}) that one has:
$$\chi(Y)\ge (b-1)(g-1)$$
with  equality holding if and only if the fibration $f$ is smooth and isotrivial. (Recall that a fibration is said to be isotrivial if all its smooth fibres are isomorphic).\newline
  If $f$ is not smooth and isotrivial, then one defines (cf. \cite{xiao}) the {\em slope} $\la(f)$ of
$f$ as:
$$\la(f):=\frac{K^2_Y- 8(b-1)(g-1)}{\chi(Y)-(b-1)(g-1)}.$$
\medskip

The fundamental inequality for the slope of a fibration is the following:

\begin{thm}[Xiao, Cornalba--Harris]\label{slope}
Let $f\colon Y\to B$ be a relatively minimal fibration, not smooth and isotrivial, with fibres of
genus $g\ge 2$. Then one has the following inequalities for the slope $\la(f)$ of $f$:
$$4(g-1)/g\le \la(f)\le 12.$$
\end{thm}
\begin{proof} See \cite{xiao}, Theorem 2. For the case of a semistable fibration, see also \cite{ch}.
\end{proof}

 We are going to reduce the proof of the Severi inequality in the general case to Theorem
\ref{slope} by means of the following:

\begin{prop}\label{trick}Let $S$ be a smooth minimal complex  surface of general type and of maximal
Albanese dimension. There exists a sequence of smooth   
complex surfaces
$Y_d$, $d\in \N_{>0}$, and  relatively minimal fibrations $f_d\colon Y_d\to \pp^1$ such that:
\begin{enumerate}
\item $f_d$ is  not isotrivial;
\item   $\lim_{d\to\infty}g_d=+\infty$, where $g_d$ is the genus of a general fibre of $f_d$;
\item $\lim_{d\to\infty}\la(f_d)=K^2_S/\chi(S)$.
\end{enumerate}
\end{prop}
\begin{proof}
Let $a\colon S\to A$ be the Albanese map and let $q\ge 2$ be the irregularity of $S$. Let $L$ be a
very ample line bundle on $A$, let $H$ be the pull-back of $L$ on $S$ and set $\al:=H^2$,
$\be:=K_SH$. 
Fix an integer $d\ge 1$,  let $\mu\colon A\to A$ be the multiplication by $d$ and
consider the following cartesian diagram:
\begin{equation}
\begin{CD}
S' @>p>>  S\\
@Va'VV @VVaV\\
A@>\mu>>A
\end{CD}
\end{equation}
The surface $S'$ is smooth minimal of maximal Albanese dimension with invariants:
\begin{equation}\label{invariants}
K^2_{S'}=d^{2q}K^2_S, \quad \chi(S')=d^{2q}\chi(S).
\end{equation}
 By \cite[Ch. 2, Prop. 3.5]{LB},  one has the following  equivalence on $A$:
$$\mu^*L\sim_{num}d^2L.$$ Hence, denoting by $H'$ the pull-back of $L$ on $S'$ via the map
$a'$,  one has: $$p^*H\sim_{num}d^2H'.$$  
Using this remark, one computes:
\begin{equation}
{H'}^2=d^{2q-4}\al,\quad H'K_{S'}=d^{2q-2}\be.
\end{equation}
Now let $D_1, D_2\in |H'|$ be  general smooth curves. Set $C_1:=D_1+D_2$ and let $C_2\in |2H'|$ be a general curve,  so that  $C_2$ is smooth and $C_1$ and $C_2$ intersect
transversely at $4d^{2q-4}\al$ points.
Let $\epsi\colon Y_d\to S'$ be the blow up of the intersection points of $C_1$ and $C_2$. 
The invariants of $Y_d$ are the following:
\begin{equation}\label{K2}
K^2_{Y_d}=K^2_{S'}-4{H'}^2=d^{2q}K^2_S-4d^{2q-4}\al 
\end{equation}
\begin{equation}\label{chi}
\chi(Y_d)=\chi(S')=d^{2q}\chi(S)
\end{equation}
The pencil spanned by $C_1$ and $C_2$ defines a fibration $f_d\colon Y_d\to\pp^1$ which is
relatively minimal by construction. Furthermore, $f_d$ is not isotrivial, since the strict transform of $C_1$ is a singular stable fibre.  The genus $g_d$ of a general fibre of $f_d$ is equal to the genus of $C_2$, namely one has:
\begin{equation}\label{genus}
g_d=1+K_{S'}H'+2{H'}^2=1+d^{2q-2}\be+2d^{2q-4}\al.
\end{equation}
 We have $\be=K_SH>0$, since $S$ is minimal of general type and $|H|$ is free, and $q>1$, since
$S$ has maximal Albanese dimension. By these remarks, statement (ii) follows from  (\ref{genus}).  

The slope of $f_d$ is given by the following expression:

\begin{equation}
\lambda(f_d)=\frac{K^2_{Y_d}+8(g_d-1)}{\chi(Y_d)+(g_d-1)}
\end{equation}
Hence statement (iii) follows by (\ref{K2}), (\ref{chi}) and  (\ref{genus}).
\end{proof}
\begin{proof}[Proof of Thm. \ref{severi}]
If $S$ is not of general type, then the inequality follows from the classification of surfaces,
hence we may assume that $S$ is of general type.

Consider the surfaces $Y_d$ of Proposition \ref{trick}.
By Theorem \ref{slope}, we have:  
\begin{equation}\label{sloped}
\la(f_d)\ge 4(g_d-1)/g_d \quad \forall d \ge 1.
\end{equation}
By Proposition \ref{trick}, taking the limit of (\ref{sloped}) for $d\to\infty$ we get:
$$K^2_S\ge 4\chi(S).$$
\end{proof}

\bigskip

\bigskip

Rita Pardini

Dipartimento di Matematica, 
Universit\`a di Pisa

Via F. Buonarroti 2,
56127 Pisa (Italy)

Email: 
pardini@dm.unipi.it
\end{document}